\documentclass[final,5p,times,twocolumn]{elsarticle}

\usepackage{amssymb}

\usepackage{amsmath}
\usepackage{amsthm}
\usepackage{bm}
\usepackage{subcaption}
\usepackage{natbib}
\usepackage{here}
\usepackage{color}
\usepackage{mathtools}
\usepackage{algorithm}
\usepackage{algpseudocode}
\mathtoolsset{showonlyrefs=true}

\theoremstyle{definition}

\newtheorem{asm}{Assumption}

\newtheorem{thm}{Theorem}
\newtheorem{lmm}{Lemma}

\newtheorem{prp}{Proposition}
\theoremstyle{remark}
\newtheorem{rmk}{Remark}

\newcommand{\e}[1]{\begin{equation}#1\end{equation}}

\newcommand{\al}[1]{\begin{align}#1\end{align}}
\newcommand{\ald}[1]{\begin{aligned}#1\end{aligned}}

\newcommand{\red}[1]{\textcolor{red}{#1}}
\newcommand{\prn}[1]{\left(#1\right)}

\newcommand{\cur}[1]{\left\{#1\right\}}
\newcommand{\ang}[1]{\left\langle#1\right\rangle}

\def\coloneq{\mathrel{\mathop:}=}
\def\r{\mathbb{R}}
\def\rn{\mathbb{R}^n}
\def\rm{\mathbb{R}^m}
\def\rd{\mathbb{R}^d}
\def\sn{\mathbb{S}^n}
\def\sd{\mathbb{S}^d}

\def\tp{\mathrm{T}}
\allowdisplaybreaks

\begin{document}

\begin{frontmatter}



\title{A feasible smoothing accelerated projected gradient method for nonsmooth convex optimization}


\author[label1]{Akatsuki Nishioka\corref{cor1}}
\cortext[cor1]{Corresponding author}
\ead{akatsuki_nishioka@mist.i.u-tokyo.ac.jp}
\affiliation[label1]{organization={Department of Mathematical Informatics, The University of Tokyo},
            addressline={Hongo 7‑3‑1}, 
            city={Bunkyo‑ku},
            postcode={113-8656}, 
            state={Tokyo},
            country={Japan}}

\author[label1,label2]{Yoshihiro Kanno}

\affiliation[label2]{organization={Mathematics and Informatics Center, The University of Tokyo},
            addressline={Hongo 7‑3‑1}, 
            city={Bunkyo‑ku},
            postcode={113-8656}, 
            state={Tokyo},
            country={Japan}}

\begin{abstract}
Smoothing accelerated gradient methods achieve faster convergence rates than that of the subgradient method for some nonsmooth convex optimization problems. However, Nesterov's extrapolation may require gradients at infeasible points, and thus they cannot be applied to some structural optimization problems. We introduce a variant of smoothing accelerated projected gradient methods where every variable is feasible. The $O(k^{-1}\log k)$ convergence rate is obtained using the Lyapunov function. We conduct a numerical experiment on the robust compliance optimization of a truss structure.
\end{abstract}




\end{frontmatter}



\section{Introduction}

Smoothing accelerated gradient methods for nonsmooth convex optimization \cite{bot15,tran17,wu23} utilize smooth approximations of nonsmooth objective functions and Nesterov's accelerated gradient method \cite{nesterov83}. Their convergence rates are $O(k^{-1}\log k)$ ($O(k^{-1})$ \cite{tran17} and $o(k^{-1}\log k)$ \cite{wu23} are achieved in some cases), which are faster than the $O(k^{-1/2}\log k)$ rate of the subgradient method \cite{beck17}. In this paper, we only consider smoothing methods that adaptively update the smoothing parameters and have global convergence guarantees. For smoothing methods with the fixed smoothing parameters, see \cite{beck17,nesterov05}.

Some problems in structural optimization have objective functions defined only on (or near) the feasible sets \cite{kanno15}. For example, when the optimization variables are some sort of sizes, it is not surprising that the objective function cannot be defined for a structure with negative variables. Smoothing accelerated projected gradient methods in \cite{bot15,tran17,wu23} are all based on Nesterov's extrapolation, which may require gradients at infeasible points. Therefore they cannot be applied to such problems.

In this short paper, we introduce a new smoothing accelerated projected gradient method where every variable is feasible. The algorithm is based on the accelerated gradient method \cite[Algorithm 20]{daspremont21} for smooth convex optimization. The $O(k^{-1}\log k)$ convergence rate is obtained using the Lyapunov function and the proof is simpler than those of \cite{bot15,tran17,wu23}. We consider the robust compliance minimization problem of a truss structure as an example, check that the problem satisfies specific assumptions, and conduct a numerical experiment.

We use the following notation: $\mathbb{R}^m_{>0}$ is the set of $m$-dimensional real vectors with positive components. $\mathbb{S}^n$, $\mathbb{S}^n_{\succeq 0}$, and $\mathbb{S}^n_{\succ 0}$ are the sets of $n$-dimensional real symmetric matrices, positive semidefinite matrices, and positive definite matrices, respectively. The inner product for $X,Y\in\sn$ is defined by $\ang{X,Y}=\mathrm{tr}(XY)$. $\|\bm{x}\|$ and $\ang{\bm{x},\bm{y}}$ denotes the Euclidean norm and the standard inner product of vectors $\bm{x},\bm{y}\in\rm$.

\section{Problem setting}

We consider a nonsmooth convex optimization problem
\e{
\underset{\bm{x}\in S}{\mathrm{Minimize}}\ \ f(\bm{x})
\label{p}
}
with the following assumptions.
\begin{asm}\label{asm1}
\ \\ \vspace{-4mm}
\begin{enumerate}
\item $S\subseteq\rm$ is a nonempty closed convex set where the projection is easily computed.
\item $f:S\to\r$ is convex and possibly nonsmooth.
\item There exists an optimal solution $\bm{x}^*\in S$.
\item A smooth approximation $f_\mu$ of $f$ satisfying the following condition is available.
\begin{enumerate}
\item For any fixed $\mu>0$, $f_\mu$ is convex and continuously differentiable on $S$.
\item There exist $L>0$, $L'\ge 0$ such that $\nabla f_\mu$ is $(L'+L/\mu)$-Lipschitz continuous on $S$.
\item There exists $\beta>0$ such that 
\e{
0 \le f_{\mu_2}(\bm{x})-f_{\mu_1}(\bm{x})\le\beta(\mu_1-\mu_2)
\label{slip}
}
for any $\bm{x}\in S$ and $\mu_1\ge\mu_2\ge0$ (we interpret $f_0=f$).
\end{enumerate}
\end{enumerate}
\end{asm}

We can construct such smoothing functions using conjugate functions, the Moreau envelope, etc. (see \cite{beck17,nesterov05,toyoda23n,tran17} for example).

\section{Motivating example}

We consider the robust (or worst-case) compliance minimization problem of a truss structure \cite{bental97,kanno15}. The optimization (or design) variable $\bm{x}\in\rm$ consists of the cross-sectional areas of bars of a truss structure, where $m$ is the number of bars. The displacement vector $\bm{u}\in\rd$ is a vector consisting of displacements of each node of a truss structure, which is computed by
\e{
K(\bm{x})\bm{u}=\bm{f}
}
where, $\bm{f}\in\rd$ is the applied force and $K(\bm{x})\in\sd_{\succeq 0}$ is the global stiffness matrix, which is a linear function ($K(\bm{x})=\sum_{j=1}^m x_j K_j$ with constant symmetric matrices $K_j\in\sd_{\succeq 0}\ (j=1,\ldots,m)$). When every component of $\bm{x}$ is positive, $K(\bm{x})$ is positive definite. By minimizing a compliance $\bm{f}^\tp\bm{u}=\bm{f}^\tp K(\bm{x})^{-1}\bm{f}$, we can obtain a stiff truss structure. Here, we want to minimize the worst-case (maximum) compliance under the uncertain load $\bm{f}=Q\hat{\bm{f}}$ ($\hat{\bm{f}}\in\rn$, $\|\hat{\bm{f}}\|=1$, $Q\in\mathbb{R}^{d\times n}$, and generally $n\le d$). By maximizing a compliance $\hat{\bm{f}}^\tp Q^\tp K(\bm{x})^{-1}Q\hat{\bm{f}}$ over $\|\hat{\bm{f}}\|=1$, we can write the robust compliance minimization problem as a minimization problem of the maximum eigenvalue \cite{kanno15,takezawa11}:
\e{\ald{
& \underset{\bm{x}\in\rm}{\mathrm{Minimize}} & & \lambda_1(Q^\tp K(\bm{x})^{-1}Q)\\
& \mathrm{subject\ to} & & \bm{l}^\tp\bm{x}\le V_0,\\
&                      & & 0<x_\mathrm{min}\le x_j\ \ (j=1,\ldots,m),
\label{p_robust}
}}
where $\lambda_1(X)$ is the maximum eigenvalue of $X$, a constant vector $\bm{l}$ consists of the lengths of bars of a truss structure, which are positive values, and $V_0>0$ is the upper limit of the volume of a truss structure. Thus, $\bm{l}^\tp\bm{x}\le V_0$ is the volume constraint (or equivalently the mass constraint with uniform density). $x_\mathrm{min}$ is the lower limit of the cross-sectional areas of bars, which ensures the invertibility of $K(\bm{x})$. This problem is a convex optimization problem since the objective function is a maximum of affine functions 
\e{
\lambda_1(Q^\tp K(\bm{x})^{-1}Q) = \underset{\|\hat{\bm{f}}\|=1}{\max}\ \underset{\bm{u}\in\rd}{\max}\  \cur{2\hat{\bm{f}}^\tp Q^\tp\bm{u}-\bm{u}^\tp K(\bm{x})\bm{u}}.
}

Obviously, the objective function of \eqref{p_robust} can be defined only on $\rm_{>0}$. The inverse of the stiffness matrix $K(\bm{x})^{-1}$ behaves completely different manner if $\bm{x}<\bm{0}$ like the function $f(x)=1/x$ and makes no physical sense.

Problem \eqref{p_robust} can also be formulated as a semidefinite programming problem, and a medium-size problem can be efficiently solved by the interior-point method \cite{bental97,kanno15}. However, for large-scale problems and possibly for other applications, we develop a fast first-order optimization method, smoothing accelerated projected gradient method with feasible variables.

Let us check that the problem \eqref{p_robust} satisfies Assumption \ref{asm1}. The first three assumptions (i)--(iii) are trivial. We consider the following smooth approximation of the maximum eigenvalue (see e.g., \cite{chen04,nesterov07}):
\e{
f_{\mu}(\bm{x})\coloneq\mu\log\left(\sum_{i=1}^n \exp{\left(\frac{\lambda_i(A(\bm{x}))}{\mu}\right)}\right)-\mu\log n
\label{approx}
}
with $A(\bm{x})=Q^\tp K(\bm{x})^{-1}Q\in\mathbb{S}^{n}$, where $\lambda_i(X)$ denotes the $i$-th largest eigenvalue of $X\in\sn$.

\begin{prp}
Let $A(\bm{x})=Q^\tp K(\bm{x})^{-1}Q\in\sn$ and $S$ be the feasible set of \eqref{p_robust}, which is compact. The smooth approximation $f_{\mu}$ defined by \eqref{approx} is convex and $\nabla f_{\mu}$ is $(L'+L/\mu)$-Lipschitz continuous for some $L,L'>0$ over $S$. Moreover, it satisfies $0 \le f_{\mu_2}(\bm{x})-f_{\mu_1}(\bm{x})\le(\log n)(\mu_1-\mu_2)$ for any $\bm{x}\in S$ and $\mu_1\ge\mu_2\ge0$.
\end{prp}
\begin{proof}
Convexity follows by the formula \cite[Example 7.16]{beck17}
\al{
f_{\mu}(\bm{x})&=\underset{X\in\sn_{\succeq 0},\mathrm{tr}(X)=1}{\sup}\cur{\ang{X,Q^\tp K(\bm{x})^{-1}Q}-\mu\sum_{i=1}^n\lambda_i(X)\log\lambda_i(X)}\\
&\ \ -\mu\log n
}
and 
\al{
\ang{X,Q^\tp K(\bm{x})^{-1}Q}
& = \sum_{i=1}^n \theta_i\bm{w}_i^\tp Q^\tp K(\bm{x})^{-1}Q\bm{w}_i\\
& = \sum_{i=1}^n \theta_i\underset{\bm{u}_i\in\rd}{\max}\  \cur{2\bm{w}_i^\tp Q^\tp\bm{u}_i-\bm{u}_i^\tp K(\bm{x})\bm{u}_i},
}
where $0\le\theta_i\le1$ and $\bm{w}_i$ are eigenvalues and orthonormal eigenvectors of $X$ satisfying $\sum_{i=1}^n \theta_i=1$. Since $f_{\mu}(\bm{x})$ is constructed by applying convexity-preserving operations (a nonnegative combination and sup operations) to affine functions, it is convex.

Next, we check the smoothness. Since the log-sum-exp function $\mu\log\left(\sum_{i=1}^n \exp{\left(x_i/\mu\right)}\right)$ is twice continuously differentiable, $\mu\log\left(\sum_{i=1}^n \exp{\left(\lambda_i(X)/\mu\right)}\right)$ is also twice continuously differentiable due to the theory of symmetric spectral functions \cite{lewis01,lewis03}. Moreover, $Q^\tp K(\bm{x})^{-1}Q$ is twice continuously differentiable on $\rm_{>0}$, and thus $f_{\mu}$ is twice continuously differentiable. We set $g(X)\coloneq \sum_{i=1}^n \exp\prn{\lambda_i(X)}\in\r$ and denote the gradient with respect to $X\in\sn$ by $\nabla_X g$. We have $\nabla_X g(X)=\sum_{i=1}^n \exp\prn{\lambda_i(X)}\bm{v}_i\bm{v}_i^\tp$, where $\bm{v}_i$ is the unit eigenvector corresponding to $\lambda_i(X)$ \cite{nesterov07}. By the similar argument to \cite{nesterov07}, for any $\bm{h}\in\rm$ we have 
\al{
& f_{\mu}(\bm{x}) = \mu\log(g(A(\bm{x})/\mu)),\\
& \ang{\nabla f_{\mu}(\bm{x}),\bm{h}}=\frac{1}{g(A(\bm{x})/\mu)}\ang{\nabla_X g(A(\bm{x})/\mu),\nabla A(\bm{x})\bm{h}},
}
and thus
\al{
& \ang{\nabla^2 f_{\mu}(\bm{x})\bm{h},\bm{h}}\\
& = -\frac{1}{\mu g(A(\bm{x})/\mu)^2}\ang{\nabla_X g(A(\bm{x})/\mu),\nabla A(\bm{x})\bm{h}}^2\\
& \ \ + \frac{1}{\mu g(A(\bm{x})/\mu)}\ang{\nabla^2_X g(A(\bm{x})/\mu)\nabla A(\bm{x})\bm{h},\nabla A(\bm{x})\bm{h}}\\
& \ \ + \frac{1}{g(A(\bm{x})/\mu)}\ang{\nabla_X g(A(\bm{x})/\mu),\bm{h}^\tp\nabla^2 A(\bm{x})\bm{h}},
\label{hessian}
}
where $\nabla A(\bm{x})\bm{h}=\sum_{j=1}^m h_j\frac{\partial A(\bm{x})}{\partial x_j}\in\sn$ and $\bm{h}^\tp\nabla^2 A(\bm{x})\bm{h}=\sum_{j,k=1}^m h_j h_k\frac{\partial^2 A(\bm{x})}{\partial x_j \partial x_k}\in\sn$. The first term of the right-hand side of \eqref{hessian} is nonpositive. By \cite[Equation (10)]{nesterov07}, the second term of \eqref{hessian} is bounded by
\e{
\frac{1}{g(A(\bm{x})/\mu)}\ang{\nabla_X g(A(\bm{x})/\mu),\bm{h}^\tp\nabla^2 A(\bm{x})\bm{h}} \le \frac{1}{\mu}\|\bm{\lambda}(\nabla A(\bm{x})\bm{h})\|_\infty^2,
}
where $\bm{\lambda}(X)=(\lambda_1(X)\ldots\lambda_n(X))^\tp$. The third term of \eqref{hessian} is bounded by
\al{
& \frac{1}{g(A(\bm{x})/\mu)}\ang{\nabla_X g(A(\bm{x})/\mu),\bm{h}^\tp\nabla^2 A(\bm{x})\bm{h}}\\
& =\frac{1}{\sum_{i=1}^n \exp\prn{\lambda_i(A(\bm{x}))/\mu}}\sum_{i=1}^n \exp\prn{\lambda_i(A(\bm{x}))/\mu}\bm{u}_i^\tp\bm{h}^\tp\nabla^2 A(\bm{x})\bm{h}\bm{u}_i\\
& \le\|\bm{\lambda}(\bm{h}^\tp\nabla^2 A(\bm{x})\bm{h})\|_\infty,
}
where $\bm{u}_i$ is the unit eigenvector corresponding to $\lambda_i(A(\bm{x}))$. Therefore, $\nabla f_{\mu}$ is $(L'+L/\mu)$-Lipschitz continuous over $S$ with
\e{
L=\underset{\bm{x}\in S}{\max}\underset{\|\bm{h}\|=1}{\max}\|\bm{\lambda}(\nabla A(\bm{x})\bm{h})\|_\infty^2,\ \ L'=\underset{\bm{x}\in S}{\max}\underset{\|\bm{h}\|=1}{\max}\|\bm{\lambda}(\bm{h}^\tp\nabla^2 A(\bm{x})\bm{h})\|_\infty,
}
both of which are finite by the compactness of $S$ and twice continuous differentiability of $A(\bm{x})$ on $S$.

Lastly, we check the relation \eqref{slip}. It follows by the argument for the log-sum-exp function \cite{toyoda23n}. We denote $\lambda_i(A(\bm{x}))$ by $\lambda_i$ and $s_i=\exp\left(\frac{\lambda_i}{\mu}\right)\ge0$ for short. We obtain
\al{
\frac{\partial}{\partial \mu}f_{\mu}(\bm{x})
& = \frac{\partial}{\partial \mu} \prn{\mu\log\left(\sum_{i=1}^n \exp\prn{\frac{\lambda_i}{\mu}}\right)}-\log n\\
& = \log\left(\sum_{i=1}^n s_i\right) + \mu\frac{\sum_{i=1}^n - \lambda_i s_i/\mu^2}{\sum_{i=1}^n s_i}-\log n\\
& = \log\prn{\sum_{i=1}^n s_i} + \frac{\sum_{i=1}^n s_i(-\log s_i)}{\sum_{i=1}^n s_i}-\log n\\
& \ge \log\prn{\sum_{i=1}^n s_i} -\log\prn{\frac{\sum_{i=1}^n s_i^2}{\sum_{i=1}^n s_i}}-\log n\\
& = \log\prn{\frac{\prn{\sum_{i=1}^n s_i}^2}{\sum_{i=1}^n s_i^2}}-\log n\\
& \ge -\log n,
}
where the fourth inequality follows by Jensen's inequality (see e.g., \cite{beck17}) for $-\log x$ with coefficients $s_i/\sum_{j=1}^n s_j$ ($i=1,\ldots,n$). Moreover,
\al{
\frac{\partial}{\partial \mu}f_\mu(\bm{x})
& = \log\prn{\sum_{i=1}^n s_i} + \frac{\sum_{i=1}^n s_i(-\log s_i)}{\sum_{i=1}^n s_i}-\log n\\
& = \frac{\sum_{i=1}^n s_i\log\prn{\sum_{j=1}^n s_j}}{\sum_{i=1}^n s_i} + \frac{\sum_{i=1}^n s_i(-\log s_i)}{\sum_{i=1}^n s_i}-\log n\\
& = \frac{1}{\sum_{i=1}^n s_i} \sum_{i=1}^n s_i \log\prn{\frac{\sum_{j=1}^n s_{j}}{s_i}} -\log n\\
& \le \log\prn{\sum_{i=1}^n \frac{s_i}{\sum_{j=1}^n s_j}\frac{\sum_{j=1}^n s_j}{s_i}} -\log n\\
& = \log n -\log n\\
& = 0,
}
where the third inequality follows by Jensen's inequality for $\log x$ with coefficients $s_i/\sum_{j=1}^n s_j$ ($i=1,\ldots,n$). Therefore, by integrating the above inequalities with respect to $\mu$ from $\mu_2$ to $\mu_1$, we obtain \eqref{slip} with $\beta=\log n$.
\end{proof}

Note that $L$ and $L'$ obtained in the above proof can be very large. Indeed, the local Lipschitz constant of $f_{\mu}$ can change drastically since the objective function behaves like the function $f(x)=1/x$. A practical (possibly nonmonotone) stepsize strategy is in the future work.

\section{Proposed method}

We consider the smoothing accelerated projected gradient (S-APG) method for problem \eqref{p} described in Algorithm \ref{sapg}, where $\Pi_S$ is the projection operator onto $S$.

\begin{algorithm}[ht]
    \caption{Smoothing accelerated projected gradient (S-APG) method with feasible variables}
    \label{sapg}
    \begin{algorithmic}[1]
    \renewcommand{\algorithmicrequire}{\textbf{function}}
    \renewcommand{\algorithmicensure}{\textbf{input:}}
    \Ensure $\bm{x}^0=\bm{z}^0\in S$, $\mu_0>0$, $a_0=0$
    \For {$k = 0,1,2,\ldots,K$}
    \State $\mu_{k} = \mu_0(k+1)^{-1}$,
    \State $L_k=L'+L/\mu_k$,
    \State $a_{k+1} = \frac{1+\sqrt{4a_{k}^2+1}}{2}$,
    \State $\bm{y}^k = \prn{1-\frac{1}{a_{k+1}}}\bm{x}^k+\frac{1}{a_{k+1}}\bm{z}^k$,
    \State $\bm{z}^{k+1} = \mathrm{\Pi}_S\prn{\bm{z}^k-\frac{a_{k+1}}{L_k}\nabla f_{\mu_k}(\bm{y}^k)}$,
    \State $\bm{x}^{k+1} = \prn{1-\frac{1}{a_{k+1}}}\bm{x}^k+\frac{1}{a_{k+1}}\bm{z}^{k+1}$.
    \EndFor
    \end{algorithmic}
\end{algorithm}

The update scheme is based on the accelerated projected gradient method \cite[Algorithm 20]{daspremont21} with strong convexity parameter $0$ for convex smooth optimization (we just replace the gradient of the objective function in \cite[Algorithm 20]{daspremont21} by the gradient of a smooth approximation). Unlike Nesterov's extrapolation $\bm{y}^{k+1}=\bm{x}^{k+1}+((a_{k}-1)/a_{k+1})(\bm{x}^{k+1}-\bm{x}^{k})$ \cite{nesterov83}, every variable of Algorithm \ref{sapg} is feasible since the formulas for $\bm{y}^k$ and $\bm{x}^{k+1}$ are convex combinations.

Let us review some important properties of the coefficient $a_{k}$ \cite{beck17,daspremont21}. The coefficient $a_{k+1}$ is a solution of the quadratic equation $x^2-x-a_{k}^2$, namely 
\e{
a_{k+1}^2-a_{k+1}-a_{k}^2=0.
\label{a_k}
}
Moreover, $0=a_0\le\ldots\le a_k$ and the following bounds hold:
\e{
a_k=\frac{1+\sqrt{4a_{k-1}^2+1}}{2}\ge\frac{1}{2}+a_{k-1}\ge\frac{1}{2}+\frac{1}{2}+a_{k-2}\ge\ldots\ge\frac{k}{2}
}
and
\e{
a_k\le\frac{1+\sqrt{4(a_{k-1}+1)^2}}{2}=a_{k-1}+\frac{3}{2}\le\ldots\le \frac{3k}{2}.
}

We write $f_{\mu_k}=f_k$ for short. We define a Lyapunov function for S-APG by
\e{
\mathcal{E}_k \coloneq \frac{a_{k}^2}{L_k}\prn{f_k(\bm{x}^{k})-f_k(\bm{x}^*)+\beta\mu_k}+\frac{1}{2}\|\bm{z}^{k}-\bm{x}^*\|^2.
}
We prove the convergence rate of S-APG using the following property of the Lyapunov function.

\begin{lmm}
For $k=0,1,\ldots$, the sequence generated by S-APG (Algorthm \ref{sapg}) for problem \eqref{p} under Assumption \ref{asm1} satisfies
\e{
\mathcal{E}_{k+1} \le \mathcal{E}_k+\frac{\beta a_{k+1}\mu_k}{L_k}.
}
\label{l_lyap}
\end{lmm}
\begin{proof}
We apply the Lyapunov argument in \cite[Theorem 4.22]{daspremont21} to $f_k$ with fixed $\mu_k$ and modify it by connecting $f_k$ and $f_{k+1}$. We have
\al{
& a_{k+1}^2\prn{f_{k}(\bm{x}^{k+1})-f_k(\bm{x}^*)} - a_{k}^2\prn{f_{k}(\bm{x}^{k})-f_k(\bm{x}^*)}\\
& = a_{k+1}\prn{f_{k}(\bm{y}^{k})-f_k(\bm{x}^*)} + a_{k}^2\prn{f_{k}(\bm{y}^{k})-f_{k}(\bm{x}^{k})}\\
& \ \ + a_{k+1}^2\prn{f_{k}(\bm{x}^{k+1})-f_{k}(\bm{y}^{k})}\\
& \le a_{k+1}\ang{\nabla f_{k}(\bm{y}^{k}),\bm{y}^k-\bm{x}^*} + a_{k}^2\ang{\nabla f_{k}(\bm{y}^{k}),\bm{y}^k-\bm{x}^k}\\
&\ \ + a_{k+1}^2\prn{\ang{\nabla f_{k}(\bm{y}^{k}),\bm{x}^{k+1}-\bm{y}^k}+\frac{L_k}{2}\|\bm{x}^{k+1}-\bm{y}^k\|^2}\\
& = \ang{\nabla f_{k}(\bm{y}^{k}),a_{k+1}^2\bm{x}^{k+1}-a_{k+1}\bm{x}^*-a_{k}^2\bm{x}^k} + \frac{a_{k+1}^2 L_k}{2}\|\bm{x}^{k+1}-\bm{y}^k\|^2\\
& = \ang{\nabla f_{k}(\bm{y}^{k}),(a_{k+1}^2-a_{k+1})\bm{x}^k+a_{k+1}\bm{z}^{k+1}-a_{k+1}\bm{x}^*-a_{k}^2\bm{x}^{k}}\\
& \ \ + \frac{L_k}{2}\|\bm{z}^{k+1}-\bm{z}^k\|^2\\
& = a_{k+1}\ang{\nabla f_{k}(\bm{y}^{k}),\bm{z}^{k+1}-\bm{x}^*} + \frac{L_k}{2}\|\bm{z}^{k+1}-\bm{z}^k\|^2\\
& \le L_k\ang{\bm{z}^{k+1}-\bm{z}^k,\bm{x}^*-\bm{z}^{k+1}} + \frac{L_k}{2}\|\bm{z}^{k+1}-\bm{z}^k\|^2\\
& = \frac{L_{k}}{2}\|\bm{z}^{k}-\bm{x}^*\|^2 - \frac{L_{k}}{2}\|\bm{z}^{k+1}-\bm{x}^*\|^2,
}
where the first relation follows by \eqref{a_k}, the second relation follows by the convexity of $f_k$ and the $L_k$-Lipschitz continuity of $\nabla f_k$ (see e.g., \cite{beck17}), the third relation follows by \eqref{a_k}, the fourth relation follows by the update rule of $\bm{x}^{k+1}$ and $\bm{y}^{k}$ in Algorithm \ref{sapg}, the fifth relation follows by \eqref{a_k}, the sixth relation follows by the update rule of $\bm{z}^{k+1}$ in Algorithm \ref{sapg} and the nonexpansiveness of the projection, the seventh relation follows by the algebraic equality $2\ang{\bm{b}-\bm{a},\bm{c}-\bm{b}}+\|\bm{a}-\bm{b}\|^2=\|\bm{a}-\bm{c}\|^2-\|\bm{b}-\bm{c}\|^2$ that holds for any $\bm{a},\bm{b},\bm{c}\in\rn$. Consequently, we obtain
\al{
& \frac{a_{k+1}^2}{L_{k}}\prn{f_{k}(\bm{x}^{k+1})-f(\bm{x}^*)}+\frac{1}{2}\|\bm{z}^{k+1}-\bm{x}^*\|^2\\
& \le \frac{a_{k}^2}{L_{k}}\prn{f_{k}(\bm{x}^{k})-f(\bm{x}^*)}+\frac{1}{2}\|\bm{z}^{k}-\bm{x}^*\|^2.
\label{ineq}
}
Moreover, since $0<\mu_{k+1}\le\mu_{k}$,
\e{
f_{k+1}(\bm{x}^{k+1})-f_{k+1}(\bm{x}^*)+\beta\mu_{k+1}\le f_{k}(\bm{x}^{k+1})-f_{k}(\bm{x}^*)+\beta\mu_k
\label{ineq2}
}
holds by \eqref{slip}. Therefore, we obtain
\al{
\mathcal{E}_{k+1} 
& = \frac{a_{k+1}^2}{L_{k+1}}\prn{f_{k+1}(\bm{x}^{k+1})-f_{k+1}(\bm{x}^*)+\beta\mu_{k+1}}+\frac{1}{2}\|\bm{z}^{k+1}-\bm{x}^*\|^2\\
& \le \frac{a_{k+1}^2}{L_{k}}\prn{f_{k}(\bm{x}^{k+1})-f_k(\bm{x}^*)+\beta\mu_{k}}+\frac{1}{2}\|\bm{z}^{k+1}-\bm{x}^*\|^2\\
& \le \frac{a_{k}^2}{L_k}\prn{f_{k}(\bm{x}^{k})-f_k(\bm{x}^*)+\frac{a_{k+1}^2}{a_{k}^2}\beta\mu_{k}}+\frac{1}{2}\|\bm{z}^{k}-\bm{x}^*\|^2\\
& = \frac{a_{k}^2}{L_k}\prn{f_{k}(\bm{x}^{k})-f_k(\bm{x}^*)+\beta\mu_{k}}+\frac{1}{2}\|\bm{z}^{k}-\bm{x}^*\|^2+\frac{\beta a_{k+1}\mu_k}{L_k}\\
& = \mathcal{E}_{k}+\frac{\beta a_{k+1}\mu_k}{L_k},
}
where the second inequality follows by $L_{k+1}\ge L_k$ and \eqref{ineq2}, the third inequality follows by \eqref{ineq}, and the fourth equality follows by \eqref{a_k}.
\end{proof}

\begin{rmk}
If we define $\tilde{\mathcal{E}}_k \coloneq \frac{a_{k}^2}{L_k}\prn{f(\bm{x}^{k})-f(\bm{x}^*)+\beta\mu_k}+\frac{1}{2}\|\bm{z}^{k}-\bm{x}^*\|^2-\sum_{l=1}^k(\beta a_{l+1}\mu_l)/L_l$, then we can obtain a monotonically nonincreasing Lyapunov function $\tilde{\mathcal{E}}_{k+1}\le\tilde{\mathcal{E}}_k$, although the argument is essencially the same.
\end{rmk}

\begin{thm}
For $k=0,1,\ldots$, the sequence generated by S-APG (Algorithm \ref{sapg}) for problem \eqref{p} under Assumption \ref{asm1} satisfies
\al{
f(\bm{x}^{k})-f(\bm{x}^*)
& =\frac{2L\|\bm{x}^{0}-\bm{x}^*\|^2+6\beta\mu_0^2\log k}{\mu_0 k}\\
& \ \ +\frac{2(L'+L/\mu_0)\prn{\|\bm{x}^{0}-\bm{x}^*\|^2+(3\beta\mu_0^2/L)\log k}}{k^2},
}
i.e., the convergence rate of S-APG is $O(k^{-1}\log k)$.
\end{thm}
\begin{proof}
By Lemma \ref{l_lyap} and $f(\bm{x}^k)-f(\bm{x}^*)\le f_k(\bm{x}^k)-f_k(\bm{x}^*)+\beta\mu_k$ from \eqref{slip}, we obtain
\al{
\frac{a_{k}^2}{L_k}\prn{f(\bm{x}^{k})-f(\bm{x}^*)}\le\mathcal{E}_{k}
& \le\ldots\le\mathcal{E}_{0}+\sum_{l=0}^{k-1}\frac{\beta a_{l+1}\mu_{l}}{L_{l}}\\
& =\frac{1}{2}\|\bm{x}^{0}-\bm{x}^*\|^2+\sum_{l=0}^{k-1}\frac{\beta a_{l+1}\mu_{l}}{L_{l}},
}
where $a_0=0$ and $\bm{z}^0=\bm{x}^0$ are used. Moreover, recalling $L_k=L'+L/\mu_k$, $a_k\le 3k/2$, and $\mu_k=\mu_0(k+1)^{-1}$, we have
\al{
\sum_{l=0}^{k-1}\frac{\beta a_{l+1}\mu_{l}}{L_{l}}
& \le \frac{\beta}{L}\sum_{l=0}^{k-1}a_{l+1}\mu_{l}^2\\
& \le\frac{3\beta\mu_0^2}{2L}\sum_{l=0}^{k-1}\frac{1}{l+1}\\
& \le\frac{3\beta\mu_0^2\log k}{2L}.
}
Therefore, recalling $a_k\ge k/2$, we obtain
\al{
f(\bm{x}^{k})-f(\bm{x}^*)
& \le\frac{L_k}{2a_{k}^2}\|\bm{x}^{0}-\bm{x}^*\|^2+\frac{L_k}{a_{k}^2}\sum_{l=0}^{k-1}\frac{\beta a_{l+1}\mu_{l}}{L_{l}}\\
& \le\frac{2L_k\|\bm{x}^{0}-\bm{x}^*\|^2}{k^2}+\frac{6\beta\mu_0^2L_k\log k}{Lk^2}\\
& =\frac{2(L'+L(k+1)/\mu_0)\|\bm{x}^{0}-\bm{x}^*\|^2}{k^2}\\
& \ \ +\frac{6\beta\mu_0^2\log k(L'+L(k+1)/\mu_0)}{Lk^2}\\
& =\frac{2L\|\bm{x}^{0}-\bm{x}^*\|^2+6\beta\mu_0^2\log k}{\mu_0 k}\\
& \ \ +\frac{2(L'+L/\mu_0)\prn{\|\bm{x}^{0}-\bm{x}^*\|^2+(3\beta\mu_0^2/L)\log k}}{k^2}.\\
\label{last}
}
\end{proof}

\begin{rmk}
When $f$ itself is continuously differentiable and $\nabla f$ is $L'$-Lipschitz continuous, by putting $L=0$ and $\beta=0$ in the first inequality of \eqref{last}, we can recover the convergence rate of the accelerated gradient method for smooth convex optimization $f(\bm{x}^{k})-f(\bm{x}^*)\le 2L'\|\bm{x}^{0}-\bm{x}^*\|^2/k^2$.
\end{rmk}

\section{Numerical result}

\begin{figure}[ht]
  \centering
  \begin{tabular}{cccc}
  \begin{minipage}[t]{0.24\hsize}
    \centering
    \includegraphics[width=2.3cm]{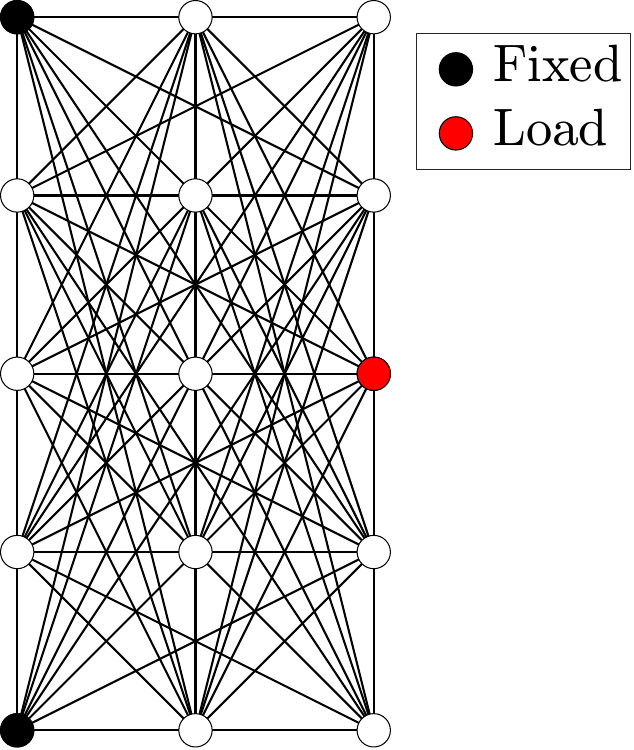}
    \vspace{-2mm}
    \subcaption{Initial design}
  \end{minipage} &
  \hspace{-4mm}
  \begin{minipage}[t]{0.24\hsize}
    \centering
    \includegraphics[width=1.4cm]{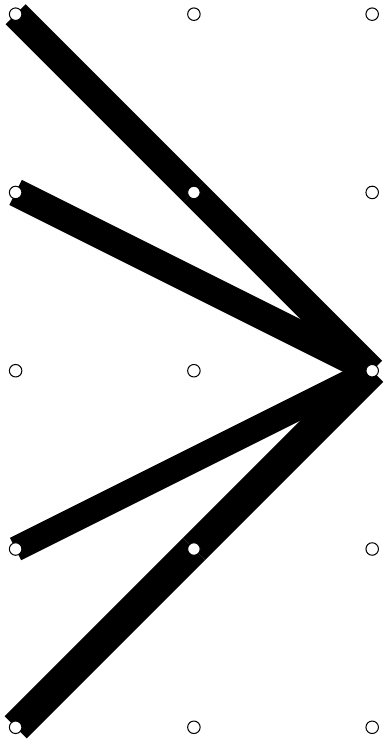}
    \vspace{-2mm}
    \subcaption{S-APG}
  \end{minipage} &
  \hspace{-5mm}
  \begin{minipage}[t]{0.24\hsize}
    \centering
    \includegraphics[width=1.4cm]{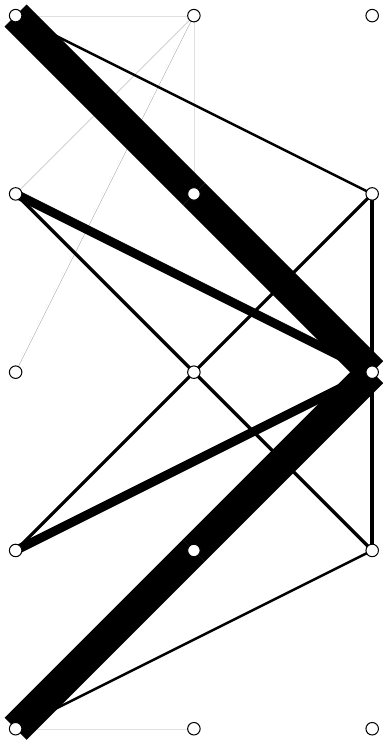}
    \vspace{-2mm}
    \subcaption{S-PG}
  \end{minipage} &
  \hspace{-5mm}
  \begin{minipage}[t]{0.24\hsize}
    \centering
    \includegraphics[width=1.4cm]{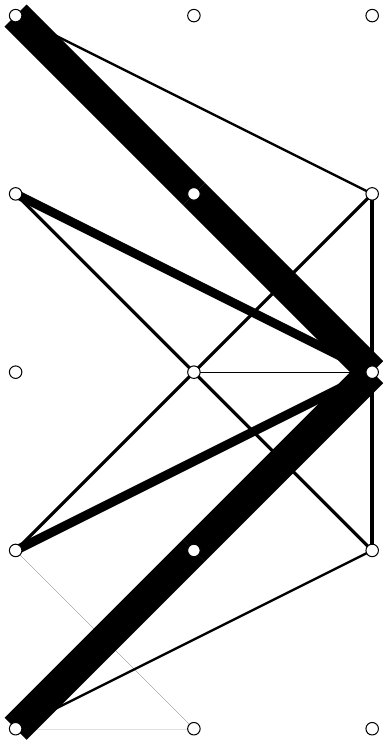}
    \vspace{-2mm}
    \subcaption{Subgrad}
  \end{minipage}
  \end{tabular}
  \caption{Comparison of the designs after 4000 iterations}
  \label{f_des}
\end{figure}

We compare S-APG with the smoothing projected gradient method without acceleration (S-PG) \cite{bian20c, toyoda23n} and the subgradient method (Subgrad) \cite{beck17} in the robust compliance minimization problem \eqref{p_robust}. We consider the initial design shown in Figure \ref{f_des}(a). The uncertain load is applied to the red node and the uncertainty set is the ellipse with the horizontal axis $2\times10^5$ N and the vertical axis $2.78\times10^5$ N. We set the parameters of the problem as follows; the number of the optimization variables is $m=74$, the size of the matrices is $n=20$, Young's modulus of the material used in the stiffness matrix is $200$ GPa, the distance between the nearest nodes is $1$ m, the upper limit of the volume is $V_0=0.1$ m$^3$, and the lower bound of the cross-sectional areas is $x_\mathrm{min}=10^{-8}$ m$^2$. We set $\mu_0=1$ and $L=10^{5}$ for S-APG and $\mu_0=1$ and $L=10^{6}$ for S-PG. The stepsize of Subgrad is $\alpha_k=10^{-6}\times k^{-1/2}$. Note that larger stepsizes lead to more oscillatory sequences in S-PG and Subgrad since stepsizes decay slower $O(k^{-1/2})$ rate compared to $O(k^{-1})$ in S-APG, and thus we set smaller stepsizes for S-PG and Subgrad than S-APG, which are experimentally fast. The experiment has been conducted on MacBook Pro (2019, 1.4 GHz Quad-Core Intel Core i5, 8 GB memory) and MATLAB R2022b.

\begin{figure}[ht]
    \centering   \includegraphics[width=6cm]{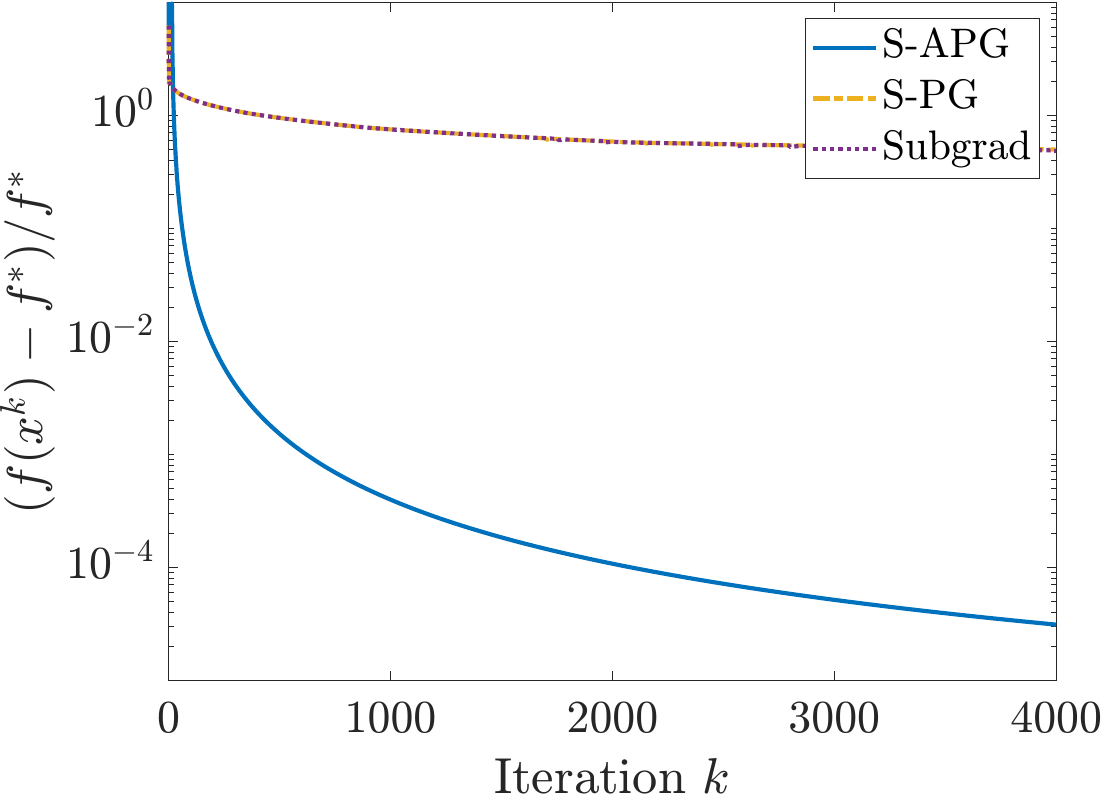}
    \caption{Difference of the objective value and the optimal value at each iteration}
    \label{f_obj}
\end{figure}

Figure \ref{f_des} shows the initial design and the designs obtained by the three algorithms after 4000 iterations (bars with cross-sectional areas less than $10^{-6}$ m$^2$ are not displayed). In this setting, the two maximum eigenvalues coincide near the optimal solution, which results in the nonsmoothness of the objective function. Figure \ref{f_obj} shows the differences between objective values and the optimal value divided by the optimal value in 4000 iterations. Note that the lines of S-PG and Subgrad overlap in Figure \ref{f_obj}. Figures \ref{f_des} and \ref{f_obj} show that S-APG converges faster than S-PG and Subgrad, both of which have $O(k^{-1/2}\log k)$ convergence rates.

\vspace{5mm}

\noindent
\textbf{Acknowledgments}

The authors are grateful to Dr.~Mirai Tanaka at The Institute of Statistical Mathematics and Dr.~Mitsuru Toyoda at Tokyo Metropolitan University for the helpful discussion. This research is part of the results of Value Exchange Engineering, a joint research project between R4D, Mercari, Inc.~and RIISE. The work of the first author is partially supported by JSPS KAKENHI JP23KJ0383. The work of the last author is partially supported by JSPS KAKENHI JP21K04351.\\






\bibliographystyle{elsarticle-harv} 
\bibliography{ref}





\end{document}